\documentclass[11pt]{article}

\usepackage{amsmath,amsfonts,amssymb} 
\usepackage[english]{babel}                        
\usepackage[active]{srcltx}                         
\usepackage{geometry}
\usepackage[T1]{fontenc}			                                

\usepackage[english]{babel}                        
\usepackage[T1]{fontenc}			   

\usepackage{theorem}{\theorembodyfont{\slshape}
  \newtheorem{theorem}{Theorem}
  \newtheorem{proposition}{Proposition}
  \newtheorem{lemma}{Lemma}
}{\theorembodyfont{\upshape}
  \newtheorem{remark}{Remark} 
}
\def\Proof{{\smallskip\noindent{\em Proof. }}}     
\def\endProof{{\hfill$\Box$\medskip\noindent}}     

\def\endProof{{\hfill$\Box$}}
\newcommand\R{{\mathbb{R}}}
\renewcommand\P{{\mathbb{P}}}
\newcommand\N{{\mathbb{N}}}

\newcommand\E{{\mathbb{E}}}

\newcommand\K{{\mathbb{K}}}

\newcommand\KK{{\mathfrak{K}}}
\newcommand\FF{{\mathfrak{F}}}
\renewcommand\S{{\mathbb{S}}}
\renewcommand\div{{\rm div}}

\renewcommand\L{{\mathbb{L}}}

\title{Fine properties of self-similar solutions of the Navier--Stokes equations\footnote{
{\it Keywords:\/} asymptotic profiles, asymptotic behavior, far-field,
large distance, selfsimilar,
incompressible viscous flows, decay estimates,
Landau stationary solutions,
homogeneous data, Oseen kernel, cancellations.\hfill\break
{\it 2000 Mathematics subject classification:\/} 76D05, 35Q30
}}

\author{
\begin{minipage}[t]{38ex}
\begin{center}Lorenzo Brandolese\end{center}\vspace*{-1.7ex}
{\sf\small%
Universit\'e de Lyon, Universit\'e Lyon~1,
CNRS UMR 5208 Institut Camille Jordan,
21 avenue Claude Bernard,
{F-69622} Villeurbanne Cedex, France.
\vspace*{-1.7ex}\begin{center}{\tt brandolese@math.univ-lyon1.fr}\end{center}}
\end{minipage}
}

\begin{document}

\maketitle

\begin{abstract}
We study the solutions of the nonstationary incompressible Navier--Stokes
equations in $\R^d$, $d\ge2$,
of  self-similar form~$u(x,t)=\frac{1}{\sqrt t}U\bigl(\frac{x}{\sqrt t}\bigr)$,
obtained from small and homogeneous initial data~$a(x)$.
We construct an explicit asymptotic formula relating the self-similar
profile~$U(x)$ of the velocity field to its corresponding initial datum $a(x)$.
\end{abstract}

\section{Introduction}

In this paper we are concerned with the study
of solutions of the elliptic problem
\begin{equation}
\label{elp}
\begin{cases}
-\frac{1}{2}U-\frac{1}{2}(x\cdot\nabla )U-\Delta U+(U\cdot\nabla U)+\nabla P=0\\
\nabla \cdot U=0,
\end{cases}
\qquad x\in\R^d,
\end{equation}
where $U=(U_1,\ldots,U_d)$ is a vector field in $\R^d$, $d\ge2$, $\nabla=(\partial_1,\ldots,\partial_d)$,
and $P$ is a scalar function defined on~$\R^d$.
Such system arises from the nonstationary Navier--Stokes equations (NS), for an incompressible viscous fluid
filling the whole $\R^d$, when looking for a velocity field  $u(x,t)$ and pressure $p(x,t)$
of forward self-similar form: $u(x,t)=\frac{1}{\sqrt t}U(x/\sqrt t)$ and  $p(x,t)=\frac{1}{t} P(x/\sqrt t)$.
An important motivation for studying the system~\eqref{elp}
is that the corresponding self-similar velocity fields $u(x,t)$
describe the asymptotic behavior at large scales for a wide class of Navier--Stokes flows.
Moreover, simple necessary and sufficient conditions for a solution
of the Navier--Stokes equations to have an asymptotically self-similar profile for large~$t$
are available, see \cite{Pla98}.
We refer to \cite{Can95} and  \cite{Lem02}, 
for more explanations and further motivations.

The problem that we address in the present paper is the study of the 
asymptotic behavior for $|x|\to\infty$ for a large class of solutions to the
system~\eqref{elp}.

\medskip
The existence of nontrivial solutions of ~\eqref{elp}
has been known for more than sixty years.
For example, in the three-dimensional case 
Landau  observed that, putting an additional axi-symmetry condition one can construct,
via ordinary differential equations methods, 
a one-parameter family  $(U,P)$, smooth outside the origin, and satisfying~\eqref{elp} in the pointwise sense for $x\not=0$
(see, e.g., \cite[p. 207]{Bat74}).

Landau's solutions have the additional property that $U$ is a homogeneous vector field
of degree~$-1$ and $P$ is homogeneous of degree~$-2$, in a such way that the
corresponding solution $(u,p)$ of (NS) turns out to be stationary.
A uniqueness result  by \v Sver\'ak, \cite{Sve06}
implies, on the other hand, that no other solution with these properties does exist
in $\R^3$, other than Landau axi-symmetric ones.
See also~\cite{KorSve} for a detailed study of the asymptotic properties of these flows.

The class of solutions to the system~\eqref{elp} is, however, much larger.
Indeed, Giga and Miyakawa \cite{GM} proposed a general method,
based on the analysis of  the vorticity equation in Morrey spaces, for constructing nonstationary
self-similar solutions of (NS).
A more direct construction was later proposed by  Cannone, Meyer, Planchon  \cite{CanMP},  \cite{CanP96},
see also~\cite[Chapt. 23]{Lem02}.
Now we know that to obtain new solutions $U$ of~\eqref{elp}
we only have to choose vector fields  $a(x)$
in  $\R^d$, homogeneous of degree $-1$, 
and satisfying some mild smallness and regularity assumption on the sphere $\S^{d-1}$:
the simplest example in $\R^3$ is obtained taking a small $\epsilon>0$
and letting
\begin{equation}
\label{simpa}
a(x)=\biggl(- \frac{\epsilon\, x_2}{|x|^2}, \frac{\epsilon\,x_1}{|x|^2}, 0\biggr),
\end{equation}
but
a condition like 
$a|_{\S^{d-1}}\in L^\infty(\S^{d-1})$ with small norm 
(or similar weaker conditions involving  the $L^d$-norm or other Besov-type
norms on the sphere) would be enough.
The basic idea is that the Cauchy problem for Navier--Stokes
can be solved,  through the application
of the contraction mapping theorem,
 in Banach spaces made of functions invariant
under the natural scaling.
The profile $U$ of the self-similar solution $u$  obtained in this way (i.e. $U=u(x,1)$)
then solves the elliptic system~\eqref{elp}.

Regularity properties and unicity classes of those (small) self-similar solutions
have been studied in different functional settings
(see, e.g.,  \cite{MiuS06}, \cite{GerPS})
and are now quite well understood.

On the other hand, probably because of 
the lack of known relations between the self-similar profile~$U$ and the datum~$a$,
even in the case of self-similar flows
emanating from the simplest data, such as in~\eqref{simpa},
the problem of the asymptotic behavior of the solutions $U$ obtained in this way
has not been addressed, before, in the literature.

The main purpose of this paper is to construct an explicit formula
relating  $U(x)$ to $a(x)$, and valid asymptotically for $|x|\to\infty$.

%

We will also consider the more general problem of 
constructing asymptotic profiles as $|x| \to\infty$
for  (not necessarily self-similar) solutions $u(x,t)$ of the Navier--Stokes equations
with slow decay at infinity
(tipically, $|u(x,t)|\le C|x|^{-1}$).
Our motivation for such generalization is that solutions with such type of decay have, in general,
a non-self-similar asymptotic for large time.
In fact, Cazenave, Dickstein and Weissler  showed that their large time 
behavior can be much more chaotic than for the solutions
described by Planchon \cite{Pla98}.
As shown in \cite{CazDW}, however, 
one can obtain some understanding on  the large time behavior
of these solutions from the analysis of their  spatial behavior at infinity.

\section{Main results and methods}
\label{sec2}

\subsection{Notations and functional spaces}

If $\mathcal{Q}=(\mathcal{Q}_{j;h,k})$ and $B=(B_{h,k})$ are, respectively, a three-order and a two-order
tensor in the Euclidean space $\R^d$, we denote by $\mathcal{Q}:B$ the vector field with components
$$(\mathcal{Q}:B)_j=
   \sum_{h,k=1}^d \mathcal{Q}_{j;h,k}B_{h,k}, \qquad j=1,\ldots,d.$$
Sometimes, in the proofs of our  decay estimates,
we will simply write $\mathcal{Q}B$ instead of $\mathcal{Q}:B$
when all components of such vectors
can be bounded by the same quantities.

We denote the gaussian function by
$$g_t(x)=(4\pi t)^{-d/2}e^{-|x|^2/(4t)}, \qquad x\in\R^d,\quad t>0.$$
As usual, we adopt the semi-group notation $e^{t\Delta}a=g_t*a$
for the solution of the heat system $\partial_t u=\Delta u$, with $u|_{t=0}=a$
for an initial datum $a$ defined on the whole~$\R^d$.

\medskip
All the functions we deal with are supposed to be measurable. 
By definition, for any $\vartheta\ge0$ and $m\in \N$,
\begin{equation}
\begin{split}
&f\in \dot E^m_\vartheta\iff
 f\in C^m(\R\backslash\{0\})\quad\hbox{and}\quad|x|^{\vartheta+|\alpha|}\, \partial^\alpha f\in L^\infty(\R^d)\quad
\hbox{$\forall\,\alpha\in\N^d,\; |\alpha|\le m$}.\\
\end{split}
\end{equation}
We are especially interested in the case $\vartheta=1$. Indeed, the spaces $\dot E^m_1$ contain
homogeneous functions of degree~$-1$
(and, in particular, the initial datum $a(x)$ given by~\eqref{simpa}).

The non-homogeneous counterpart of $\dot E^m_\vartheta$ is
the smaller space $E_\vartheta^m$, which is defined by the additional requirement that $\partial^\alpha f\in L^\infty(\R^d)$
for all $|\alpha|\le m$.
These spaces are equipped with their natural norm:
\begin{equation*}
\begin{split}
&\|f\|_{\dot E_\vartheta^m}=\max_{|\alpha|\le m}\,\sup_{x\in\R^d\backslash\{0\}} |x|^{\vartheta+|\alpha|}|\partial^\alpha f(x)|,\\
&\|f\|_{E_\vartheta^m}=\max_{|\alpha|\le m}\,\sup_{x\in\R^d} (1+|x|)^{\vartheta+|\alpha|}|\partial^\alpha f(x)|.
\end{split}
\end{equation*}

Our starting point is a classical result by Cannone, Meyer and Planchon about
the construction of self-similar solutions of the Navier--Stokes equations
\begin{equation*}
\begin{cases}
 \partial_t u+\nabla\cdot(u\otimes u)=\Delta u -\nabla p\\
 \hbox{div}\,u=0\\
 u|_{t=0}=a,
\end{cases}
\qquad x\in\R^d,\, t>0.
\end{equation*}
Even though their construction goes through under very general assumptions
on the regularity of the initial data, here we are mainly interested in the following
simple result:

\begin{theorem}[see \cite{CanMP}, \cite{CanP96}]
\label{theoCMP}
For all $m\in\N$ there exist $\epsilon,\beta>0$ such that for all initial datum $a\in \dot E^m_1$
homogeneous
of degree~$-1$, divergence-free and satisfying
\begin{equation}
\label{smallCMP}
\|a\|_{\dot E_1^m} <\epsilon,
\end{equation}
there exists a unique self-similar solution $u(x,t)=\frac{1}{\sqrt t}U\Bigl(\frac{x}{\sqrt t}\Bigr)$
of the Navier--Stokes system  (written in the usual integral form, see  (NS) below) starting from~$a$,
and such that $\|U\|_{E^m_1}<\beta$.
Moreover,
\begin{equation}
 \label{proCMP}
U(x)=e^{\Delta}a+{\cal O}(|x|^{-2}), \qquad \hbox{as $|x|\to\infty$}.
\end{equation}
\end{theorem}

More precisely, Cannone, Meyer and Planchon prove that $U(x)=e^{\Delta}a(x)+{\cal R}(x)$,
where the remainder term satisfies ${\cal R}\in E_2^m$.
Their result was stated in dimension three, but their proof easily adapts for all $d\ge2$.

\subsection{Main results}

Our main result shows that one can give a much more precise
asymptotic formula between the asymptotic profile $U(x)$ and the datum $a(x)$.
It turns out that such asymptotic profile has a different structure in different space dimensions.

\begin{theorem}
\label{theoss}
Let $a(x)$ be a homogeneous datum of degree $-1$, such that $a$ is smooth
on the unit sphere $\S^{d-1}$ and
satisfying the smallness condition~\eqref{smallCMP} for some $m\ge3$.
Let $u(x,t)=\frac{1}{\sqrt t}U\Bigl(\frac{x}{\sqrt t}\Bigr)$ the self-similar solution constructed
in Theorem~\ref{theoCMP}.
Then the following profiles hold:
\begin{subequations}
\begin{itemize}
\item If $d=2$, we have as $|x|\to\infty$,
 \begin{equation}
\label{profss2D}
U(x)=a(x)-\log(|x|)\frac{\mathcal{Q}(x)\!:\!A}{|x|^6} + \mathcal{O}(|x|^{-3}), 
\end{equation}
Here $A=(A_{h,k})$ is the $2\times 2$ matrix given by $A_{h,k}=\int_{\S^1} (a_ha_k)$ and
$\mathcal{Q}(x)=\bigl(\mathcal{Q}_{j;h,k}(x)\bigr)$, where 
the $\mathcal Q_{j;h,k}$
are  homogeneous polynomials of degree three
(given by the explicit formula~\eqref{explicitQ} below)
\item
For $d=3$, we have as $|x|\to\infty$,
\begin{equation}
\label{profss3D}
U(x)=a(x)+\Delta a(x)-\P \nabla \cdot(a\otimes a)
-\frac{ \mathcal{Q}(x)\!:\!B}{|x|^7} + \mathcal{O}\bigl(|x|^{-5}\log|x|\bigr),
\end{equation}
for a $d\times d$ constant real matrix  $B=(B_{h,k})$ depending on~$a$.
Here $\P=Id-\nabla(\Delta)^{-1}\hbox{div}$ is the Leray-Hopf projector onto the divergence-free vector fields.
\item
For $d\ge4$, the far-field asymptotics reads, as $|x|\to\infty$,
\begin{equation}
\label{profss4D}
U(x)=a(x)+\Delta a(x)-\P \nabla\cdot (a\otimes a)+ \mathcal{O}\bigl(|x|^{-5}\log|x|\bigr).
\end{equation}
\end{itemize}
\end{subequations}
\end{theorem}

In  Section~\ref{sec9} we will restate and prove this theorem
in a more general form, removing the assumption that $a$ is homogeneous.
Such more general theorem will apply also for solutions
$u(x,t)$ of Navier--Stokes of non-self-similar form.
On the other hand, we will not seek for the greatest generality about the regularity of the
datum: even though there is a considerable interest in studying
self-similar solutions emanating from rough data (see  \cite{Lem02}, \cite{Gru06}),
in most of our statements we will assume that
$a\in C^3(\S^{d-1})$, which is of course non-optimal, but permits us to greatly simplify
the presentation of our results and to better emphasize the main ideas.

The method that we present in this paper would allow to compute, in principle,
the asymptotics of $U$ up to any order,  when $a$ is smooth on $\S^{d-1}$.
However, the higher order terms have quite complicated expressions.

The functions $\Delta a$ and $\P\nabla\cdot (a\otimes a)$ 
appearing in our expansions
are both homogeneous of degree $-3$
and smooth outside the origin.
Therefore, our asymptotic profiles imply that 
conclusion~\eqref{proCMP} of Theorem~\ref{theoCMP}
can be improved into
$$U(x)=a(x)+\mathcal{O}\bigl(|x|^{-3}\log|x|\bigr), \qquad\hbox{as $|x|\to\infty$}, \qquad\hbox{if $d=2$},$$
and
$$U(x)=a(x)+\mathcal{O}\bigl(|x|^{-3}\bigr), \qquad\hbox{as $|x|\to\infty$}, \qquad\hbox{if $d\ge3$}.$$
The datum $a$ can be replaced here by its filtered version $e^{\Delta}a$.

It turns out that such improved estimates are optimal for generic 
self-similar solutions.
For example, in the two dimensional case, the logarithmic factor cannot
be removed, since the improved bound $U(x)=a(x)+\mathcal{O}\bigl(|x|^{-3}\bigr)$
would require $\mathcal{Q}\!:\!A\equiv0$: such stringent condition can be proved to be equivalent
to the orthogonality relations $\int_{\S^1} a_1^2=\int_{\S^1} a_2^2$ and $\int_{\S^1} a_1a_2=0$.

\subsection{Main methods}

We will use the  semigroup method and the theory of mild solutions of the Navier--Stokes
equations  as
explained in detail in the books \cite{Can95} and~\cite{Lem02}.
The main novelty of our approach relies on the use of the following ingredients:
\begin{enumerate}
 \item 
The first one is the use of remarkable, but not so much known,
\emph{cancellation properties\/} hidden inside the kernel $\K(x,t)$ of the Oseen operator $e^{t\Delta}\P$,
and inside other related operators, appearing in the integral formulation of Navier--Stokes.

To be more precise,  we can write $\K(x,t)=\KK(x)+t^{-d/2}\K_2(x/\sqrt t)$,
where $\KK(x)$ is a tensor whose components are homogeneous function of degree~$-d$
(namely, second order derivatives of the fundamental solution of the Laplacian in~$\R^d$),
and $\K_2$ is exponentially  decaying  as $|x|\to\infty$.
Such decomposition already played an important role in our previous work~\cite{BraV07},
where we showed that solutions $u(x,t)$ \emph{arising from well-localized data\/}
behave like 
$$u(x,t)\sim \nabla_x\KK(x):{E}(t), \qquad\hbox{as $|x|\to\infty$},$$
where ${E}(t)$ is the energy matrix of the flow: ${E}(t)=\bigl(\int_0^t\!\int u_hu_k(y,s)\,dy\,ds\bigr)$.

A crucial fact in the proof of the results  of the present paper 
will be the use  of the identities, for $j=1,\ldots,d$,
\begin{equation*}
\int_{\S^{d-1}}\KK(\omega)\,d\omega=0, 
\qquad
\int_{\S^{d-1}}\omega_j\nabla\KK(\omega)\,d\omega=0. 
\end{equation*}
Such cancellations 
are somehow hidden in $\K$, because the non-homogeneous part $\K_2$
(and, {\it a fortiori\/}, the kernel $\K$)
\emph{does not have\/} a vanishing integral on the sphere.

\item
Our second ingredient are \emph{asymptotic formulae for convolution integrals\/}:
roughly speaking, these formulae 
consist in deducing
the exact profile as $|x|\to\infty$ of a convolution
product $f*g(x)$, from information on the regularity, the cancellations, and the 
behavior at infinity of the two factors~$f$ and~$g$.
In their simplest form, and for $f$ and $g$ ``well behaved'' at infinity, those formulae read
\begin{equation}
\label{afci}
f*g(x)\sim\Bigl(\int f\Bigr)g(x)+\Bigl(\int g\Bigr)f(x),
\qquad \hbox{as $|x|\to\infty$}.
\end{equation}
We will apply several generalizations and variants of~\eqref{afci} in different situations
(including the case of non-integrable functions)
the factors $f$ and $g$ being either
the Oseen kernel, the heat kernel, or a function related to the non-linearity.
The assumptions for the validity of~\eqref{afci} are quite stringent
(notice that~\eqref{afci} is obviously wrong if, {\it e.g.\/},  $f$ and $g$ are both a gaussian function).
Neverthless, the method that we use here has a wide applicability and can be used
for constructing the far-field asymptotics for equations of other equations.
See, {\it e.g.\/}, \cite{BraK07} for an application to a class of convection equations with anomalous diffusion.
\end{enumerate}

We will also make use 
of the so called
{\it bi-integral formula\/}.
Such formula is obtained by simply iterating the usual integral formulation of the Navier--Stokes equations,
which we now recall:
\begin{equation*}
\tag\hbox{NS}
\begin{cases}
u(t)=e^{t\Delta}a-\displaystyle\int_0^t e^{(t-s)\Delta}\P\nabla\cdot (u\otimes u)(s)\,ds\\
\div(a)=0.
\end{cases}
\end{equation*}
Using the Oseen kernel $\K(x, t)$, we can define
the Navier--Stokes bilinear operator
as
\begin{equation*}
\label{rep1}
B(u,v)(t)=\int_0^t \K(t-s)*\nabla \cdot(u\otimes v)(s)\,ds.
\end{equation*}
Then (NS)  can be written simply as
$u=e^{t\Delta}a-B(u,u)$.
The bi-integral formula
is obtained by a straightforward  iteration:
\begin{equation}
\label{bif}
u(t)=e^{t\Delta}a-B(e^{t\Delta}a,e^{t\Delta}a)+2B(e^{t\Delta}a,B(u,u))
-B(B(u,u),B(u,u)).
\end{equation}

\medskip

Roughly speaking, combining equation~\eqref{bif} with some nice properties of the heat kernel and
fine decay estimates
of the bilinear operator
we can prove ({\it e.g.\/} when $d=3$) that
$e^{t\Delta}a\sim a+t\Delta a$ and $B(e^{t\Delta}a,e^{t\Delta}a)(x)\sim t\P\nabla\cdot(a\otimes a)(x)$ as
$|x|\to\infty$.
After obtaining an explicit far-field asymptotics for $u(x,t)$, it is easy to deduce, in the self-similar case,
the behavior at infinity of the profile $U(x)$, by passing to self-similar variables and eliminating~$t$.
The two last terms in the above bi-integral formula will contribute to the remaining terms in the right-hand side
of expansion~\eqref{profss3D}.

Notice that, in the two-dimensional case, the term
$\P\nabla\cdot(a\otimes a)$ is not well-defined when $a$ is homogeneous of degree~$-1$.
This explains the different structure of our asymptotic expansions in this case.
The special structure of the asymptotic profiles in the two-dimensional case
can be observed also if, instead of considering the behavior for $|x|\to\infty$ as we do in this paper,
one focuses on the the behavior of solutions for large time.
(See, {\it e.g.\/}, \cite{GW3}).

\bigskip
For sake of simplicity, in this paper we consider only data
such that $a(x)\sim |x|^{-1}$ as $|x|\to\infty$, which is the natural assumption
for the study of global strong solutions and related self-similarity phenomena.

However, the study of the asymptotic behavior for large~$|x|$ of solutions $u$
(possibly defined only locally in time) is also of interest
in more general situations, such as $a(x)\sim |x|^{-\vartheta}$.
The far-field behavior of the solution $u(x,t)$ of (NS), 
then mainly depends on the competition between three factors.
The first one is the spatial localization of the datum (say, the value of the exponent~$\vartheta$)
and the consequent space-time decay of the linear evolution $e^{t\Delta}a$.
The other two factors are the action of the quadratic non-linearity $u\otimes u$ and of the
 non-local operator $\P\hbox{div}(\cdot)$.

When $\vartheta>d+1$,
the action of this nonlocal operator  (whose kernel behaves at infinity
like $|x|^{-d-1}$) is predominant, and is responsible of \emph{spatial spreading effects}.
When $(d+1)/2<\vartheta<d+1$ (the limit case $\vartheta=(d+1)/2$
corresponding to the situations in which $u\otimes u$ decays like the kernel),
the linear evolution becomes predominant
and the spatial spreading phenomenon is not directly observed on the solution,
but rather on its~\emph{fluctuation} $u-e^{t\Delta}a$.
We refer to our previous paper~\cite{BraV07} for a sharp description of these issues.

The asymptotic profiles of $u$ as $|x|\to \infty$ in the cases $0<\vartheta<1$ and $1<\vartheta<(d+1)/2$ 
should have a slightly different structure, but they are not known with precision
yet. The method that we use in this paper for $\vartheta=1$,
and in particular the idea of iterating the Duhamel formula making use of the cancellations of the kernels,
might be used to compute them.
More and more iterations would be needed to deal with data decaying slower than $|x|^{-1}$,
or to determine the asymptotics to a higher order.
On the other hand, no iteration or cancellation property was needed for the faster decaying data studied in~\cite{BraV07}.

\bigskip
The plan of the paper is the following: we begin with the study of the Oseen kernel.
In Section~\ref{sec4}, after some generalities about the asymptotic of convolutions, we describe the behavior
at large distances of solutions to the heat equation.
Section~\ref{sec5} is devoted to the (more or less standard) construction of solutions with a prescribed
space-time decay.
In Section~\ref{section6} we show how to use the cancellations of the Oseen kernel
to get some new fine estimates.
In the remaining part of the paper we will state and prove a
more general form of Theorem~\ref{theoss}.

\section{Asymptotics and cancellations of the Oseen kernel $\K$ and of the kernel $F$}
\label{sec3}

Let $\K(x,t)$ be the kernel of $e^{t\Delta}\P$, let $F(x,t)$ be the kernel of $e^{t\Delta}\P\hbox{div}(\cdot)$.
Both $\K(\cdot,t)$ and $F(\cdot, t)$ belong to $C^\infty(\R^d)$ and they satisfy the scaling properties 
$\K(x,t)=t^{-d/2}\K(x/\sqrt t,1)$ and $F(x,t)=t^{-(d+1)/2}F(x/\sqrt t,1)$.

Denote by $\Gamma$ the Euler Gamma function and by $\delta_{j,k}$ the Kronecker symbol.
The following Proposition extends and completes a Lemma contained in~\cite{BraV07}.
\begin{proposition}
\label{asG}
Let $\KK=(\KK_{j,k})$, where $\KK_{j,k}(x)$ is the homogeneous function of degree~$-d$
\begin{subequations}
\begin{equation}
\label{KK}
\KK_{j,k}(x)=\frac{\Gamma(d/2)}{2\pi^{d/2}}\cdot \frac{\bigl(-\delta_{j,k}|x|^2+dx_jx_k\bigr)} {|x|^{d+2}},
\end{equation}
and $\FF=(\FF_{j;h,k})$, where $\FF_{j;h,k}=\partial_h\KK_{j,k}$, which we can write also as
\begin{equation}
\label{FF}
\FF_{j;h,k}(x)=\frac{\Gamma\bigl(\frac{d+2}{2}\bigr)}{\pi^{d/2}}\cdot
   \frac{\sigma_{j,h,k}(x)|x|^2-(d+2)x_jx_hx_k}{|x|^{d+4}},
\end{equation}
where $\sigma_{j,h,k}(x)=\delta_{j,h}x_k+\delta_{h,k}x_j+\delta_{k,j}x_h$, for $j,h,k=1,\ldots,d$.
\end{subequations}
\begin{subequations}
Then the following decompositions hold:
\begin{equation}
\label{prof KK}
\K(x,t)=\KK(x)+|x|^{-d}\Psi\Bigl(x/\sqrt t\Bigr),
\end{equation}
and
\begin{equation}
\label{profKF}
F(x,t)=\FF(x)+|x|^{-d-1}\widetilde\Psi\Bigl(x/\sqrt t\Bigr),
\end{equation}
where $\Psi$ and $\widetilde\Psi$ are smooth outside the origin and such that, for all $\alpha\in\N^d$, and $x\not=0$, 
$|\partial^\alpha \Psi(x)|+|\partial^\alpha \widetilde\Psi(x)|\le Ce^{-c|x|^2}$.
Here $C$ and $c$ are positive constant,
depending on $|\alpha|$ but not on $x$.
\end{subequations}

Moreover, the following cancellations hold:
\begin{equation}
\label{canc}
 \left\{
\begin{aligned}
& \displaystyle\int_{\S^{d-1}} \KK(\omega)\,d\sigma(\omega)=\displaystyle\int_{\S^{d-1}}
 \omega_\ell\KK(\omega)\,d\sigma(\omega)=0\\
& \displaystyle\int_{\S^{d-1}}\FF(\omega)\,d\sigma(\omega)
   =\displaystyle\int_{\S^{d-1}} \omega_{\ell}\FF(\omega)d\sigma(\omega)=0\\
&\displaystyle\int_{\S^{d-1}} \omega_{\ell}\omega_{m}\FF(\omega)d\sigma(\omega)=0,
\qquad\qquad\ell,m=1,\ldots,d.
\end{aligned}
\right.
\end{equation}
\end{proposition}

\begin{remark}
 \label{remQ}
The homogeneous polynomials $\mathcal{Q}(x)=\bigl(\mathcal Q_{j;h,k}(x)\bigr)$ appearing in the
statement of Theorem~\ref{theoss} is defined by the relation
$\FF(x)=|x|^{-d-4}\mathcal{Q}(x)$, that is, with  
$\gamma_d=\Gamma(\frac{d+2}{2})/\pi^{d/2}$,

\begin{equation}
 \label{explicitQ}
\mathcal{Q}_{j;h,k}(x)=\gamma_d\Bigl( (\delta_{j,h}x_k+\delta_{h,k} x_j+\delta_{k,j}x_h)|x|^2-(d+2)x_j x_h x_k\Bigr).
\end{equation}

\end{remark}

\Proof
The symbol of $\K$ is
\begin{equation*}
\label{symb}
\widehat \K_{j,k}(\xi,t)=e^{-t|\xi|^2}\biggl(\delta_{j,k}-\frac{\xi_j\xi_k}{|\xi|^2}\biggr)
=e^{-t|\xi|^2}\delta_{j,k}-\int_t^\infty \xi_j\xi_ke^{-s|\xi|^2}\,ds
\end{equation*}

Taking the inverse Fourier transform we get
\begin{equation*}
\K_{j,k}(x,t)= \delta_{j,k}\,g_t(x)+\int_t^\infty\partial_j\partial_k g_s(x)\,ds\equiv
\K^{(1)}_{j,k}(x,t)+\K^{(2)}_{j,k}(x,t) 
\end{equation*}
Computing the derivatives $\partial_j\partial_k g_s(x)$ and changing the
variable $\lambda=\frac{|x|}{\sqrt{4s}}$ in the integral we get

\begin{equation*}
\K^{(2)}_{j,k}=\pi^{-d/2}|x|^{-d}\int_0^{|x|/\sqrt {4t}}
\biggl( 
-\delta_{j,k}\, \lambda^{d-1}
+2\lambda^{d+1}\frac{x_j x_k}{|x|^2}
\biggr)
 e^{-\lambda^2}\,d\lambda.
\end{equation*}

But, for all $r>0$ and $\alpha>-1$, 
$$
\int_0^r \lambda^\alpha e^{-\lambda^2}\,d\lambda
=\frac{1}{2}\Gamma\biggl(\frac{\alpha+1}{2}\biggr)
-\int_r^\infty\lambda^\alpha e^{-\lambda^2}\,d\lambda.
$$
Choosing first $\alpha=d-1$, then $\alpha=d+1$ and using 
$\Gamma((d+2)/2)=(d/2)\Gamma(d/2)$,
we get
\begin{equation*}
\label{bG}
\K^{(2)}_{j,k}(x,t) =\frac{|x|^{-d}}{2\pi^{d/2}}\Gamma\Bigl(\frac{d}{2}\Bigr)
\biggl[-\delta_{j,k}+d\frac{x_jx_k}{|x|^2}\biggr]+
    |x|^{-d}\Psi_{j,k}(x/\sqrt {4t}).
\end{equation*}
Here, $\Psi=(\Psi_{j,k})$ is a family of  functions such that,
\begin{equation}
\label{est Psi}
\forall\,\alpha\in\N^d, \quad|\partial^\alpha \Psi(y)|\le C_\alpha e^{-c|y|^2}, \qquad y\in\R^3.
\end{equation}
Observing that $\K^{(1)}_{j,k}$ can be bounded by the second term on the right hand side
and modifying, if necessary, the functions $\Psi_{j,k}$ 
(which can be done without  affecting estimate~\eqref{est Psi}) we see that decomposition~\eqref{prof KK} holds.
The decomposition~\eqref{profKF} is now an immediate consequence of the definition of $\FF(x)$.

Observe that $\KK_{j,k}=\partial_j\partial_k \E_d$, where $\E_d$ is the fundamental solution of $-\Delta$ in~$\R^d$.
From the radial symmetry of $\E_d$, we immediately get
\begin{equation*}
\label{cancp}
 \displaystyle\int_{\S^{d-1}} \KK_{j,k}(\omega)\,d\sigma(\omega)=\displaystyle\int_{\S^{d-1}}
 \omega_j\KK(\omega)\,d\sigma(\omega)=0, \qquad j\not=k
\end{equation*}
and
\begin{equation*}
 \displaystyle\int_{\S^{d-1}}\FF(\omega)\,d\sigma(\omega)=
\displaystyle\int_{\S^{d-1}} \omega_{\ell}\omega_{m}\FF(\omega)d\sigma(\omega)=0,
\qquad\qquad\ell,m=1,\ldots,d.
\end{equation*}
Using again the radiality of $\E_d$ and 
$\Delta \E_d=0$ on $\S^{d-1}$,
yields $\int_{\S^{d-1}}\KK_{j,j}(\omega)\,d\sigma(\omega)=0$.
This argument also shows that
the identities
\begin{equation*}
\int_{\S^{d-1}}\omega_\ell\,\FF_{j;h,k}(\omega)\,d\sigma(\omega)=0,
\qquad j,h,k,\ell=1,\ldots,d
\end{equation*}
can be reduced to the proof of the equality
\begin{equation}
\label{dirre}
\int_{\S^{d-1}} \omega_\ell\,\partial_\ell\partial_j^2 \E_d(\omega)\,d\sigma(\omega)
=\int_{\S^{d-1}}\omega_\ell\, \partial_\ell^3 \E_d(\omega)\,d\sigma(\omega),
\qquad j\not=\ell.
\end{equation}
The fact that both terms in~\eqref{dirre} are zero follows from
$\partial_j\partial_h\partial_k \E_d(\omega)=\mathcal{Q}_{j,h,k}(\omega)$, for
$\omega\in\S^{d-1}$ and formula~\eqref{explicitQ}.
In the computation, one needs to use the moment relation
$$\int_{\S^{d-1}}\omega_j^2\,d\sigma(\omega)=\frac{1}{d}\int_{\S^{d-1}}\,d\sigma(\omega)$$
and the well known identities (easily obtained via the Stokes formula)
\begin{equation*}
\left\{
\begin{aligned}
&\displaystyle\int_{\S^{d-1}}\omega_j^4\,d\sigma(\omega)=\frac{3}{d(d+2)}\displaystyle\int_{\S^{d-1}}d\sigma(\omega)\\
&\displaystyle\int_{\S^{d-1}}\omega_j^2\omega_k^2\,d\sigma(\omega)=\frac{1}{d(d+2)}\displaystyle\int_{\S^{d-1}}d\sigma(\omega),
\qquad j\not=k.
\end{aligned}
\right.
\end{equation*}

\endProof

\section{Far-field asymptotics of convolutions and application to the heat equation}
\label{sec4}

The purpose of our next result is to describe the exact behavior as
$|x|\to \infty$ 
of the convolution product of two functions $f$ and $g$
from the asymptotic properties of each factor.
We will consider only a simple particular situation that will be sufficient for our purposes.

\begin{proposition}
\label{pro1}
Let $d\ge1$ and $m\ge1$ two integers.
Let $f\in \dot E^m_\vartheta$ for some $0\le \vartheta<d$
and $g\in L^1(\R^d,(1+|x|)^m dx)\cap \dot E^0_{d+m}$.
Then the convolution product $f*g$ satisfies
\begin{equation}
\label{cpp}
f*g(x)=\sum_{ \genfrac{}{}{0pt}{}{\gamma\in\N^d}{0\le|\gamma|\le m-1} } \frac{(-1)^{|\gamma|} }{\gamma!}
\biggl(\int y^\gamma g(y)\,dy\biggr)\partial^\gamma f(x) \,+\, {\cal R}(x),
\end{equation}
where  ${\cal R}(x)$ is a remainder  term satisfying,
for some constant $C>0$ independent of $f$ and $g$ and all $x\not=0$:
\begin{equation}
\label{reme}
\bigl| \mathcal{R}(x)\bigr|\le C|x|^{-m-\vartheta} \|f\|_{\dot E^m_\vartheta}
	\Bigl( \|g\|_{\dot E^0_{d+m}}+ \|g\|_{L^1(\R^d,|x|^m\,dx)}\Bigr).
\end{equation}
\end{proposition}

\begin{remark}
The identity~\eqref{cpp} is useful, for large~$|x|$, when at least one derivative $|\partial ^{\gamma}f|$
decays at infinity exactly as $c_{\gamma}|x|^{-\vartheta-\gamma}$
(at least in some directions).
In this case, ${\cal R}(x)$ is indeed a lower order term as $|x|\to\infty$.
\end{remark}

\Proof
We can assume, without restriction, that
$ \|f\|_{\dot E^m_\vartheta}= \|g\|_{\dot E^0_{m+d}}=1$.
We have to estimate the difference between $\int f(x-y)g(y)\,dy$ and the first term
on the right-hand side of~\eqref{cpp}.
Such difference can be written as the sum of four terms $D_1+\cdots+D_4$, where
\begin{equation*}
\begin{split}
&D_1\equiv\int_{|y|\le |x|/2} \Bigl[f(x-y)-\sum_{|\gamma|\le m-1} 
\frac{(-1)^{|\gamma|}}{\gamma!}\partial^\gamma f(x)y^\gamma\Bigr]g(y)\,dy,\\
&D_2\equiv\int_{|y|\le |x|/2} g(x-y)  f(y)\,dy,\\
&D_3\equiv\int_{|y|\ge |x|/2,\; |x-y|\ge |x|/2} f(x-y)g(y)\,dy\,dy
\end{split}
\end{equation*}
and
$$D_4\equiv - \sum_{|\gamma|\le m-1} 
\frac{(-1)^{|\gamma|}}{\gamma!}\partial^\gamma f(x)
\int_{|y|\ge|x|/2} y^\gamma g(y)\,dy.$$

Using the Taylor formula, we see that
$$ |D_1|\le C |x|^{-\vartheta-m}\int_{|y|\le |x|/2}|y|^m|g(y)|\,dy,$$
which is bounded by the right-hand side of~\eqref{reme}.
Direct estimates show that
$|D_2|$, $|D_3|$ and $|D_4|$ are bounded by  $C|x|^{-\vartheta-m}$ as well.

\endProof

\medskip
\begin{remark}
We can give now
a more precise statement about the asymptotics claimed in~\eqref{afci}.
The simplest result reads as follow: if $f,g\in E^1_{\alpha+d}$ (the non-homogeneous space)
for some $\alpha>0$, then
$$f*g(x)=\Bigl(\int f\Bigr)g(x)+\Bigl(\int g\Bigr)f(x)+\mathcal{O}\bigl(|x|^{-d-\alpha^*}\bigr),$$
 as  $|x|\to\infty$,
where $\alpha^*=\min\{2\alpha,\alpha+1\}$.
When $\alpha=1$ the remainder must be replaced by $\mathcal{O}\big(|x|^{-d-2}\log(|x|)\bigr)$.
The proof relies on the same argument as that used in the proof of Proposition~\ref{pro1}:
the only difference is that the Taylor formula is applied to both $f$ and $g$, so that one has to introduce
an additional term $D_5$ in the decomposition of $f*g$.
Of course, one could state many variants of this result: here the most important
condition was  that the decay of the two factors $f$ and $g$ (or at least the decay of 
the factor with the least spatial localization)
must increase after derivation; but one could put, instead, a more general condition in terms
of moduli of continuity.
\end{remark}

Other  useful functional spaces are, for $m\in\N$, $\vartheta\ge0$,
\begin{equation*}
X^m_\vartheta =
\Bigl\{ u\in L^1_{\rm loc}((0,\infty), C^m(\R^d))\colon \;  \|u\|_{X_\vartheta^m}
\equiv \max_{|\alpha|\le m}\,
\hbox{ess}\sup_{\!\!\!\!\!\!\!\!\! x,t}\,(\sqrt t+|x|)^{\vartheta+|\alpha|}|\partial^\alpha_x u(x,t)|<\infty \Bigr\}.
\end{equation*}

The use of such spaces for the Navier--Stokes equations
is more or less classical (see, {\it e.g.\/}, \cite{CanP96}, \cite{CazDW})
but, unfortunately, there is no agreement on the notations.

\bigskip
\medskip

The following lemma is elementary:

\begin{lemma}
\label{heat l}
Let $m\in\N$, $a\in\dot E^m_\vartheta$, with $0\le\vartheta<d$.
Then there is a constant $C>0$, independent on $a$, such that
\begin{equation*}
\label{heat est}
\|e^{t\Delta}a\|_{X^m_\vartheta}\le C\|a\|_{\dot E^m_\vartheta}.
\end{equation*}
\end{lemma}

\Proof
For $\alpha\in\N^d$, $|\alpha|\le m$, one writes 
$\partial_x^\alpha e^{t\Delta}a(x)=\int \partial_x^\alpha g_t(x-y)\,a(y)\,dy$
and splits the integral in $\R^d$ into the three new integrals,
corresponding to the three disjoint regions $|y|\le|x|/2$, $|x-y|< |x|/2$ and the complementary region in $\R^d$.
For the second integral one first applies $|\alpha|$-times integration by parts. Then the direct estimate
$|\partial_x^\alpha g_t(x)|\le C|x|^{-d-|\alpha|}$
 gives the spatial decay $|\partial_x^\alpha e^{t\Delta}a(x)|\le C|x|^{-\vartheta-|\alpha|}$.
On the other hand, $a$ belongs to the Lorentz space $L^{d/\vartheta,\infty}(\R^d)$,
and $\partial_x^\alpha g_t\in L^{d/(d-\vartheta),1}(\R^d)$.
Then the time decay estimate $\|\partial_x^\alpha e^{t\Delta}a\|_\infty\le Ct^{-(\vartheta+|\alpha|)/2}$
follows from the generalized Young inequality (see, {\it e.g.\/} \cite{Lem02}).

\endProof

\medskip
As an application, we get the exact asymptotic profile as $|x|\to\infty$ for the solution of the Cauchy
problem associated with the heat equation for slowly oscillating data.
We first recall two standard notations: if $\beta\in\N^d$ we set:
$(2\beta-1)!!=\prod_{ \genfrac{}{}{0pt}{}{j=1,\ldots,d}{\beta_j\ge1} } 1\cdot3\cdot\ldots\cdot(2\beta_j-1)$ and
$(2\beta)!!=\prod_{ \genfrac{}{}{0pt}{}{j=1,\ldots,d}{\beta_j\ge1} } 2\cdot 4\cdot \ldots\cdot 2\beta$.
Now we can state the following:
\begin{lemma}
\label {as heat lem}
(i)
Let $m\ge1$ be an integer, $0\le\vartheta<d$ and $a\in \dot E^m_\vartheta$.
Then,
\begin{equation*}
 \label{hke}
e^{t\Delta}a(x)=\sum_{|2\beta|\le m-1} \frac{(2t)^\beta}{(2\beta)!!} \,\partial^{2\beta}a(x)\,
+\,\mathcal{O}\bigl(t^{m/2}|x|^{-\vartheta-m}\bigr), \qquad\hbox{as $|x|\to\infty$,}
\end{equation*}
uniformly for $t>0$
({\it i.e.\/} the remainder term is bounded by $Ct^{m/2}|x|^{-\vartheta-m}$).

(ii) 
In particular, if $m\ge4$ and $a\in \dot E^m_1$:
\begin{equation*}
 \label{as heat lap}
e^{t\Delta}a(x) =a(x)+t\Delta a(x) +\mathcal{O}\bigl(t^2|x|^{-5}\bigr), \qquad
\hbox{as $|x|\to\infty$},
\end{equation*}
uniformly for $t>0$.
\end{lemma}

\Proof
Indeed, writing $e^{t\Delta}a(x)=g_t*a(x)$, we can apply Proposition~\ref{pro1} 
with $g_t(x)=(4\pi t)^{-d/2}e^{-|x|^2/(4t)}$ instead of $g$.
Observing that,  for all $\beta\in\N^d$,
$$\int y^{2\beta}g_t(y)\,dy=2^\beta(2\beta-1)!!\,t^\beta,$$
and that $\int y^\gamma g_t(y)\,dy=0$ if $\gamma\in\N^d$ is not of the form $\gamma=2\beta$,
we obtain the result.

\endProof

\section{Global existence of decaying solutions}
\label{sec5}

We already recalled that $F$ denotes the kernel of $e^{t\Delta}\P\hbox{div}(\cdot)$
and, for $t>0$,
$F(x,t)=t^{-(d+1)/2}F(x/\sqrt t,1)$.
It is also well known that $|\partial_x^\alpha F(x,1)|\le C_\alpha(1+|x|)^{-d-1-|\alpha|}$ for all $\alpha\in\N^d$.
A quick way to prove this decay  at infinity is to observe that such estimate
is immediate for $|x|\ge1$ for both terms in the right-hand side of equation~\eqref{profKF}.
Moreover, it is clear from its definition that $F(\cdot,t)\in C^\infty(\R^d)$ for $t>0$.

\medskip
Let us introduce the linear operator $\L$, defined on $d\times d$ matrices $w=(w_{h,k})$ by the relation
\begin{equation}
\label{Lw}
\L(w)(x,t)=\int_0^t\!\!\int F(x-y,t-s)w(y,s)\,dy\,ds.
\end{equation}
More explicitly (and accordingly with the notation introduced in Section~\ref{sec2}),
 the $j$-component is given by
\begin{equation*}
\label{Lwj}
\L(w)_j(x,t)=\int_0^t\!\!\int \sum_{h.k} F_{j;h,k}(x-y,t-s)w_{h,k}(y,s)\,dy\,ds.
\end{equation*}

The interest of considering such operator is that the Navier--Stokes bilinear operator can be expressed as
\begin{equation*}
 \label{thi}
B(u,v)=\L(u\otimes v).
\end{equation*}

We start with a simple lemma (already known in a slightly less general
form, see \cite{Miy02}, \cite{CazDW}).

\begin{lemma}
 \label{sll}
Let $m\in\N$ and $w=(w_{h,k})\in X^m_2$.
Then $\L(w)\in X^m_1$ and, for some constant $C>0$ independent of $w$,
\begin{equation}
\label{Lcont}
\|\L(w)\|_{X^m_1}\le C\|w\|_{X^m_2}.
\end{equation}
\end{lemma}

\Proof
We can assume, with no loss of generality, $\|w\|_{X^m_2}=1$.
We start writing
\begin{equation}
\label{dB}
\partial^\alpha_x\L(w)(t)=\int_0^t\!\!\int \partial^\alpha_x
F(x-y,t-s)w(y,s)\,ds.
\end{equation}
Let $\alpha\in\N^d$, such that $|\alpha|\le m$ and $x\not=0$.
We split the spatial integral in equation~\eqref{dB} into the three regions
$|y|\le |x|/2$, next $|x-y|\le |x|/2$, then ($|y|\ge|x|/2$ and $|x-y|\ge|x|/2$)
and we denote with~$I_1$, $I_2$ and~$I_3$, the three corresponding integrals.
From the estimate (deduced from~\eqref{profKF}) 
 $|\partial^\alpha_x F(x,t)|\le C|x|^{-(d+|\alpha|)}t^{-1/2}$,
and the estimate $|w(y,s)|\le |y|^{-1}s^{-1/2}$,
we obtain immediately
\begin{equation}
\label{twe}
|I_1(x,t)|+|I_3(x,t)|\le C|x|^{-1-|\alpha|}.
\end{equation}
We now treat $I_2$.
When~$|\alpha|=0$ we can simply use well known fact that $\|F(\cdot,t)\|_1\le Ct^{-1/2}$
to obtain $|I_2(x,t)|\le C|x|^{-1}$. When $1\le |\alpha|\le m$
we make as many integration by parts as needed,
and use estimates of the form (deduced from the  rescaling
properties of $F$ recalled at the beginning of this section and the fact that $\partial^\alpha_x F(\cdot,1)\in L^1(\R^d,(1+|x|)^{|\alpha|}\,dx)$\,)
$$\|\, |\!\cdot\!|^\alpha \partial^\alpha_x F(\cdot,t-s)\|_1 \le C(t-s)^{-1/2}.$$
Then observing that $|\partial_y^\alpha w(y,s)|\le |y|^{-1-|\alpha|}s^{-1/2}$ for $|\alpha|\le m$,
we conclude that $I_2$ can be estimated like $I_1$ and $I_3$ in~\eqref{twe}.
Summarizing,
we showed that
\begin{equation}
\label{spd}
\bigl|\partial^\alpha_x \L(w)(x,t) \bigr|\le C|x|^{-1-|\alpha|}.
\end{equation}

There is a now well known strategy (see \cite{Miy02})
to deduce time decay estimates from the corresponding
space decay estimates. 
Namely, using the semi-group property of the Oseen kernel,
$$ \L(w)(t)=e^{t\Delta/2}\L(w)(t/2)+\int_{t/2}^t F(t-s)*w(s)\,ds\equiv K_1(t)+K_2(t).$$
From the Young inequality in Lorentz spaces,
and observing that $\bigl\|  \L(w)(t)\bigr\|_{L^{d,\infty}}$
is uniformly bounded, because of inequality~\eqref{spd}, we get
$$
\bigl\|\partial^\alpha_x K_1(t)\bigr\|_\infty
\le \bigr\|\partial^\alpha_x g_{t/2}\bigr\|_{L^{d/(d-1),1}} \bigl\| \L(w)(t/2)\bigr\|_{L^{d,\infty}}
\le C\,t^{-(1+|\alpha|)/2}.
$$
Moreover,
\begin{equation*}
\begin{split}
\bigl\|\partial_x^\alpha K_2(t)\bigr\|_\infty
\le \int_{t/2}^t \|F(t-s)\|_1\|\partial_x^\alpha w(s)\|_\infty\,ds\le C\,t^{-(1+|\alpha|)/2}.
\end{split}
\end{equation*}

Concluding, we showed that
$$\bigl |\partial^\alpha_x \L(w)\bigr|(x,t)
\le C\bigl( |x|^{-(1+|\alpha|)}\wedge t^{-(1+|\alpha|)/2}\bigr)\le C' (\sqrt t+|x|)^{-1-|\alpha|}.$$
This proves the natural estimate~\eqref{Lcont}.

\endProof

\bigskip

We now follow the standard procedure for constructing
global solutions to (NS) in the space $X^m_1$.
Our starting point will be the following basic existence result,
which is nothing but a reformulation of well-know results in the literature 
(see \cite{Can95}, \cite{CanMP}, \cite{CazDW}, \cite{Miy02})
in a slightly more general form.

\begin{proposition}
\label{theorem1}
Let $d\ge2$ and  $m\ge0$ be two integers.
There exist two constants $\epsilon>0$ and $M>0$ such that for all divergence-free vector field
$a\in \dot E^m_1$, satisfying
\begin{equation*}
\|a\|_{\dot E^m_1} <\epsilon,
\end{equation*}
there exists a unique solution $u\in X^m_1$ of (NS) starting from $a$
(in the sense that $u(t)\to a$ in $\mathcal{S}'(\R^d)$, as $t\to0$),
such that $\|u\|_{X^m_1}\le \epsilon M$. 
\end{proposition}

\Proof
We only  have to apply the size estimate for the linear evolution
$$\|e^{t\Delta} a\|_{X^m_1}\le C\|a\|_{\dot E^m_1}$$
(this is a particular case of Lemma~\ref{heat l})
and the corresponding estimate for the bilinear operator:
\begin{equation*}
\label{bicont}
\|B(u,v)\|_{X^m_1}\le C\|u\|_{X^m_1} \|v\|_{X^m_1}.
\end{equation*}
This last inequality is obtained applying Lemma~\ref{sll} with $w=u\otimes v$.
The existence a solution~$u\in X^m_1$ (and its unicity in a ball of such space)
now follows from the application of the contraction mapping theorem,
as explained {\it e.g.\/} in Cannone's book~\cite{Can95}.
Slightly changing the estimates of the previous Lemma we easily obtain, {\it e.g.\/},
the bound $|B(u,u)(x,t)|\le C|x|^{-3/2}t^{1/4}$, implying $B(u,u)(t)\to0$ in~$\mathcal{S}'(\R^d)$
as $t\to0$. Thus, from (NS), $u(t)\to a$ as $t\to0$ in the distributional sense.

\endProof

\begin{remark}
In the particular case in which~$a$ is a homogeneous vector field
of degree $-1$ in~$\R^d$,
the solution~$u$ constructed in~Proposition~\ref{theorem1} is self-similar:
\begin{equation*}
u(x,t)=\frac{1}{\sqrt t} U\biggl(\frac{x}{\sqrt t}\biggr),
\end{equation*}
for some with $U\in E^m_1$ (the non-homogeneous space).
This easily follows from the scaling invariance of (NS)
(see {\it e.g.\/}~\cite[Ch. 3]{Can95}).
\end{remark}

\section{Fine estimates of the bilinear term}
\label{section6}

It follows from Lemma~\ref{sll} that, for $w\in X^0_2$, we have
\begin{equation}
\label{trin}
|\L(w)(x,t)|\le C(t^{-1/2}\wedge |x|^{-1}).
\end{equation}
This was enough for constructing a decaying solution of (NS).

However, to obtain such decay estimate we used only few properties of the kernel $F(x,t)$,
namely, its pointwise decay and its rescaling properties.
Next Lemma will allow us to considerably improve estimate~\eqref{trin}, at least in the
parabolic region $|x|\ge \sqrt t$.
Its proof will make an essential use of the  \emph{cancellations properties\/} of the kernel $F(x,t)$
and requires some regularity for~$w$.

\begin{lemma}
\label{fine estim}
Let
$w=(w_{h,k})$, with $w\in X^2_2$.
Let $\L(w)$ be defined by equality~\eqref{Lw}.
\begin{subequations}
Then  we have, for $d\ge3$,
\begin{equation}
 \label{decay Lw}
|\L(w)(x,t)|\le C\Bigl(t^{-1/2}\wedge \, t\,|x|^{-3}\Bigr).
\end{equation}
When $d=2$, we have the weaker estimate
\begin{equation}
 \label{fine est2}
|\L(w)(x,t)|\le Ct|x|^{-3}\log\Bigl(\frac{|x|}{\sqrt t}\Bigr), \qquad |x|\ge e\sqrt t.
\end{equation}
\end{subequations}
Under the more stringent assumption $w\in X^3_2$, we have the following estimates for $\nabla\L(w)$:
\begin{equation*}
\label{grad Lw}
|\nabla \L(w)(x,t)|\le 
\begin{cases}
C\bigl(t^{-1}\wedge t|x|^{-4}\bigr), &\hbox{if $d\ge3$}\\
Ct^{-1} &\hbox{if $d=2$ and $|x|\le e\sqrt t$}\\
Ct|x|^{-4}\log\bigl(|x|/\sqrt t\bigr) &\hbox{if $d=2$ and $|x|\ge e\sqrt t$}.
\end{cases}
\end{equation*}
In all these inequalities $C>0$ is a constant dependent on $w$
only through its $\|\cdot\|_{X^2_2}$ or its $\|\cdot\|_{X^3_2}$-norm,
and independent on~$x$ and~$t$.
\end{lemma}

\Proof
We can limit ourselves to the region $|x|\ge e\sqrt t$.
Indeed, when $|x|\le e\sqrt t$ the result holds because of inequality~\eqref{Lcont},
which, in the special case $m=0,1$, implies
$|\L(w)(x,t)|\le Ct^{-1/2}$ and $|\nabla \L(x)(x,t)|\le Ct^{-1}$.

Let us decompose
\begin{equation}
\begin{split}
\label{dec L} 
\L(w)(x,t)&=\int_0^t\!\!\int_{|y|\le |x|/2}F(x-y,t-s)w(y,s)\,dy\,ds \\
		&\qquad+\int_0^t\!\!\int_{|y|\le |x|/2}F(y,t-s)w(x-y,s)\,dy\,ds\\
		&\qquad+\int_0^t\!\!\int_{|y|\ge |x|/2,\;|x-y|\ge |x|/2}F(x-y,t-s)w(y,s)\,dy\,ds\\
		&\equiv\L_1+\L_2+\L_3
\end{split}
\end{equation}
We start with estimating $\L_3$.
Using $|F(x-y,t-s)|\le C|x-y|^{-d-1}\le C'|y|^{-d-1}$ (the two inequalities being valid in the region of~$\R^d$
where we perform the integration) and $|w(y,s)|\le |y|^{-2}$,
we get $|\L_3(x,t)|\le Ct|x|^{-3}$.

In view of the use of the Taylor formula, we further decompose $\L_1$ 
(recalling also~\eqref{profKF})
as
\begin{equation}
\label{L1 dec}
\begin{split}
	\L_1=& \int_0^t\!\!\int_{|y|\le |x|/2}\bigl[F(x-y,t-s)-F(x,t-s)\bigr]w(y,s)\,dy\,ds\\
   &\qquad + \FF(x)\!:\!\int_0^t \int_{|y|\le |x|/2} w(y,s)\,dy\,ds\\
  &\qquad +|x|^{-d-1}\int_0^t \widetilde \Psi(x/\sqrt{t-s})\int_{|y|\le |x|/2} w(y,s)\,dy\,ds.
\end{split}
\end{equation}
Using $|\nabla F(x,t)|\le C|x|^{-d-2}$, next $|y|\,|w(y,s)|\le C|y|^{-1}$ 
shows that the first term in~\eqref{L1 dec} is bounded
by $Ct|x|^{-3}$.

When $d\ge3$, since $|w(y,s)|\le C|y|^{-2}$,
the second term in the right-hand side of~\eqref{L1 dec} is also bounded by $Ct|x|^{-3}$.
When $d=2$, make use of the inequality $|w(y,s)|\le C(\sqrt s+|y|)^{-2}$ and of the chage
of variables $y=\sqrt s z$.
This leads to the weaker upper bound estimate
of the form $Ct|x|^{-3}\log(|x|/\sqrt t)$, valid for $|x|\ge e\sqrt t$.

The simplest way to treat the third term on the right-hand side of~\eqref{L1 dec}
is to recall that  $|\widetilde \Psi (x)|\le C$.
In this way, one can proceed exactly as for the previous term and obtain the same
bounds.
This would be enough for the proof of this Lemma.
However, for later use
(namely, to shorten the proof of Lemma~\ref{lle} below), we want to prove that this last term in~\eqref{L1 dec}
is bounded, in the region $|x|\ge e\sqrt t$,  by $Ct|x|^{-3}$ also when $d=2$. This is easy: indeed
$\widetilde\Psi$ has a fast decay at infinity; here, the use of the inequality $|\Psi(x)|\le C|x|^{-1}$
is enough to conclude.

We now consider $\L_2$. We decompose it as
\begin{equation}
\label{L2 dec}
\begin{split}
\L_2=  &\int_0^t\!\!\int_{|y|\le |x|/2}  F(y,t-s) \bigl[w(x-y,s)-w(x,s)+y\cdot\nabla w(x,s)\bigr]\,dy\,ds\\
  &\qquad\int_0^t w(x,s)\int_{|y|\le |x|/2} F(y,t-s)\,dy\,ds \\
  &\qquad - \int_0^t \nabla w(x,s)\cdot \int_{|y|\le |x|/2}y \,F(y,t-s)\,dy\,ds.\\
\end{split}
\end{equation}
Now we use the inequalities $|\nabla^2 w(x,t)|\le C|x|^{-4}$ 
and $|y|^2\,|F(y,t-s)|\le C|y|^{-d+1}$,
and obtain that the first term on the right-hand side
in~\eqref{L2 dec} is bounded by $Ct|x|^{-3}$.
We now conclude using the cancellations of the kernel~$F$:
more precisely, since $\int F(\cdot,t-s)\,dy=0$ and $|F(y,t-s)|\le C|y|^{-d-1}$ the second term
is also bounded by $Ct|x|^{-3}$.

A brutal estimate of the third term in~\eqref{L2 dec} would give a non-optimal bound
of the form $C|x|^{-3}\log(|x|\sqrt t)$ for large $|x|$, which is not enough.
But, for $|x|\ge 2\sqrt{t}$, the third term in~\eqref{L2 dec}
can be further decomposed as
\begin{equation}
\label{L23 dec}
\begin{split}
&\int_0^t \nabla w(x,s)\cdot \int_{|y|\le \sqrt{t-s}} y\,F(y,t-s)\,dy\,ds\\
 & \qquad+\int_0^t \nabla w(x,s)\cdot \int_{\sqrt{t-s}\le |y|\le |x|/2} y\,\FF(y)\,dy\,ds\\
 & \qquad+\int_0^t \nabla w(x,s)\cdot \int_{\sqrt{t-s}\le |y|\le |x|/2} y\,|y|^{-d-1}\widetilde \Psi(y/\sqrt{t-s})\,dy\,ds.
 \end{split}
\end{equation}
Now it is easy to see that the first and the third term in~\eqref{L23 dec} are $\mathcal{O}(t|x|^{-3})$.
But $\FF$ has vanishing first order moments on the sphere (see Proposition~\ref{asG})
so that the second term in~\eqref{L23 dec} is zero.

Summarizing, we have established inequality~\eqref{fine est2} in the two-dimensional case and
inequality~\eqref{decay Lw} when $d\ge3$.

\medskip
To prove the inequality for~$\nabla\L$,
we fix $\ell\in\{1,\ldots,d\}$ and we write
\begin{equation*}
 \partial_\ell L(x,t)=\int_0^t\!\! F(x-y,t-s)\partial_\ell w(y,s)\,dy\,ds\equiv \widetilde \L_1+\widetilde\L_2+\widetilde\L_3,
\end{equation*}
where the decoposition is obtained as before (see~\eqref{dec L}).
The two terms $\widetilde \L_2$ and $\widetilde \L_3$ are treated exactly as before, but
we get now upper bound of the form $Ct|x|^{-4}$ since $\partial_\ell w$ (and its derivatives up to the second order)
decays faster than $w$ (and its corresponding derivatives).
Notice that we need use here the assumption $w\in X^3_2$ which ensures a decay for the derivatives
up to the order three.

For treating~$\widetilde \L_1$ we integrate by parts.
It is easy to see that the boundary term is bounded by $Ct|x|^{-4}$.
The other term is $\int_0^t\int_{|y|\le |x|/2}\partial_\ell F(x-y,t-s)w(y,s)\,dy\,ds$,
for which we obtain the usual bound $Ct|x|^{-4}$ when $d\ge3$ and
$Ct|x|^{-4}\log(|x|/\sqrt t)$ for $d=2$ and $|x|\ge e\sqrt t$.

\endProof

\medskip
\begin{remark}
 \label{rem fine}
For later use, let us observe that if
$u\in X^2_1$ is the solution constructed in Proposition~\ref{theorem1}, in the case $m\ge2$,
then, applying Lemma~\ref{fine estim} to $w=u\otimes u$, so that
$\L(w)=B(u,u)$, we get
\begin{subequations}
\begin{equation}
 \label{buud}
|B(u,u)|(x,t)\le
\begin{cases}   
 C(t^{-1/2}\wedge t|x|^{-3}) &\qquad\hbox{if  $d\ge3$}\\
 Ct^{-1/2} &\qquad\hbox{if $d=2$ and $|x|\le e\sqrt t$}\\
Ct|x|^{-3}\log(|x|/\sqrt t) &\qquad \hbox{if $d=2$ and $|x|\ge e\sqrt t$}.
\end{cases}
\end{equation}
In the case $u\in X^3_1$ (this requires the more stringent assumption $a\in \dot E^3_1$
in Proposition~\ref{theorem1}), in addition to the above estimates, the bilinear term satisfies
\begin{equation}
\label{grad Buu}
|\nabla B(u,u)|(x,t)\le 
\begin{cases}
C(t^{-1}\wedge t|x|^{-4}), &\hbox{if $d\ge3$}\\
Ct^{-1} &\hbox{if $d=2$ and $|x|\le e\sqrt t$}\\
Ct|x|^{-4}\log(|x|/\sqrt t) &\hbox{if $d=2$ and $|x|\ge e\sqrt t$}.
\end{cases}
\end{equation}
\end{subequations}
These estimates will play an essential role in the study of the bi-integral formula
\end{remark}
\begin{equation}
\label{biff}
u(t)=e^{t\Delta}a-B(e^{t\Delta}a,e^{t\Delta}a)+2B(e^{t\Delta}a,B(u,u))
-B(B(u,u),B(u,u)).
\end{equation}


\section{Asymptotic profiles of  the velocity field in the 2D case}

In the two-dimensional case, 
from Lemma~\ref{heat l} and Remark~\ref{rem fine} we get, for $(x,t)\in\R^2\times(0,\infty)$,
\begin{equation}
 \label{qw}
|e^{t\Delta}a\otimes B(u,u)(x,t)|\le
 \begin{cases}
Ct^{-1} & \hbox{if $|x|\le e\sqrt t$}\\
Ct|x|^{-4}\log(|x|/\sqrt t) &\hbox{if $|x|\ge e\sqrt t$}.
\end{cases}
\end{equation}
The last term in~\eqref{biff} satisfies, always for $(x,t)\in\R^2\times(0,\infty)$, an even stronger estimate, namely
\begin{equation}
 \label{qwB}
|B(u,u)\otimes B(u,u)(x,t)|\le
 \begin{cases}
Ct^{-1} & \hbox{if $|x|\le e\sqrt t$}\\
Ct^2|x|^{-6}\log^2(|x|/\sqrt t) &\hbox{if $|x|\ge e\sqrt t$}.
\end{cases}
\end{equation}

Next Lemma allows us to show that the two last terms in the right-hand side
of~\eqref{biff} can be considered as remainders, {\it i.e.\/}, they can be included in the $\mathcal{O}(t|x|^{-3})$ term.

\begin{lemma}
\label{le5}
Let $w=(w_{h,k})$ defined on $\R^2\times (0,\infty)$
with $w_{h,k}(x,t)$ bounded by the right hand side of~\eqref{qw},
or by the right-hand side of~\eqref{qwB}.
Then, if~$\L(w)$ is given by~\eqref{Lw}, we have for some $C>0$ independent on $x$ or~$t$,
$$|\L(x,t)|\le C\Bigl(t^{-{1/2}}\wedge \,t\,|x|^{-3}\Bigr).$$
\end{lemma}

\Proof
Our assumptions imply $w\in X^0_2$.
Then we deduce from Lemma~\ref{sll} that $|\L(x,t)|\le Ct^{-1/2}$, therefore we can assume that 
$|x|\ge e\sqrt t$.
Then we split the spatial integral defining $\L$ (see~\eqref{Lw}) into the three regions
$|y|\le \sqrt s$, $\sqrt s\le |y|\le |x|/2$ and $|y|\ge |x|/2$.
The first term that we obtain is bounded using $|F(x-y,t-s)|\le C|x|^{-3}$ (this is true only in 2D)
and $|w(y,s)|\le Cs^{-1}$. For the second term we use the same bound for $F$ and $|w(y,s)|\le Cs|y|^{-4}\log(|y|/\sqrt s)$.
The last term is treated  using the bound $|w(y,s)|\le C\sqrt s|y|^{-3}$ and 
that $\|F(t-s)\|_1\le C(t-s)^{-1/2}$.

\endProof

\medskip
Next Lemma will be useful for treating the term $B(e^{t\Delta}a,e^{t\Delta}a)$ arising in~\eqref{biff}.
Note that for $a\in \dot E^2_1$ we have, from Lemma~\ref{heat l}, $e^{t\Delta}a \otimes e^{t\Delta}a\in X^2_2$.

\begin{lemma}
\label{lle}
Let $w=(w_{h,k})$, with $w_{h,k}\in X^2_2$ for all $h,k=1,2$.
Then we have
\begin{equation}
\label{exp Lw}
\L(w)(x,t)=\FF(x)\!:\!\int_0^t\!\!\int_{|y|\le|x|} w(y,s)\,dy\,ds+\mathcal{O}(t|x|^{-3}),
\qquad\hbox{as $|x|\to\infty$},
\end{equation}
uniformly with respect to $t$ in the region $|x|\ge e \sqrt t$.
Here $\L(w)$ is given by~\eqref{Lw} and $\FF(x)$ is the homomeneous tensor 
of order three defined by
equation~\eqref{FF}.
\end{lemma}

\Proof
This follows from the proof of Lemma~\ref{fine estim}.
Therein, we decomposed $\L(w)$ as the sum of several terms, all of which, excepted one,
could be bounded by
$Ct|x|^{-3}$.
The only term for which such upper bound could brake down was
\begin{equation*}
 \FF(x)\!:\!\int_0^t\int_{|y|\le|x|/2} w(y,s)\,dy\,ds
\end{equation*}
(see the second term in the right-hand side of~\eqref{L1 dec}).
A simple modification of the error term now shows that we can
change the above domain of the spatial integral into $\{|y|\le |x|\}$.

\endProof

\begin{lemma}
\label{lol}
Let $a(x)$ be a vector field defined on $\R^2$, such that $a\in \dot E^2_1$.
Then, for $|x|\to\infty$ and uniformly in time, in the region  $|x|\ge e\sqrt t$, we have:
\begin{equation*}
\label{inth2}
\int_0^t\!\!\int_{|y|\le |x|} (e^{s\Delta}a\otimes  e^{s\Delta}a)(y)\,dy\,ds
=\int_0^t\int_{\sqrt s\le |y|\le |x|} (a\otimes a) (y)\,dy\,ds+\mathcal{O}(t\,1)
\end{equation*}
(here and below $\mathcal{O}(t\,1)$ denotes a remainder function bounded by $Ct$ for $|x|\ge e\sqrt t$).

In particular, if $a$ is homogeneous in $\R^2$ of degree~$-1$
\begin{equation*}
\label{int hh}
\int_0^t\!\!\int_{|y|\le |x|} (e^{s\Delta}a\otimes e^{s\Delta}a)(y)\,dy\,ds
= t\log\Bigl(\frac{|x|}{\sqrt t}\Bigr) \biggl(\int_{\S^1} a\otimes a\biggr)+\mathcal{O}(t\,1),
\quad\hbox{ as $|x|\to\infty$},
\end{equation*}
\end{lemma}

\Proof
Indeed, we can assume $|x|\ge \sqrt t$.
Then $\displaystyle\int_0^t\!\!\displaystyle\int_{|y|\le \sqrt s} (e^{s\Delta}a\otimes e^{s\Delta}a)(y)\,dy\,ds$ is bounded
by $Ct$.
It remains to treat  
$$\int_0^t\!\!\int_{\sqrt s\le |y|\le |x|} (e^{s\Delta}a\otimes e^{s\Delta}a)(y)\,dy\,ds,$$
which we can rewrite as the sum of four new integrals,
if we use the decomposition $e^{t\Delta}a(x)=a(x)+\mathcal{R}(x,t)$ obtained in Lemma~\ref{as heat lem}
(in the case $\vartheta=1$, $m=2$)
and a similar decomposition for $e^{t\Delta}b$.
Here, $\mathcal{R}$ satisfies
$|\mathcal{R}(x,t)|\le Ct|x|^{-3}$.
An easy calculation shows that the three integrals containing at least one factor $\mathcal{R}$ are bounded
by $Ct$.

\endProof

\begin{theorem}
\label{theo2D}
Let $u(x,t)\in X^2_1$ be the global solution of the Navier--Stokes equations in $\R^2$, with datum $a\in \dot E^2_1$
(as constructed in Proposition~\ref{theorem1}).
Then $u$ has the following profile for $|x|\to\infty$, uniformly with respect to~$t$ in the region~$|x|\ge e\sqrt t$:
\begin{equation}
 \label{pro2Dnss}
u(x,t)=a(x)-\FF(x)\!:\!\int_0^t\!\!\int_{\sqrt s\le |y|\le |x|} (a\otimes a)(y)\,dy\,ds\;+\; \mathcal{O}(t|x|^{-3}).
\end{equation}
Moreover, if $a$ is homogeneous of degree $-1$, then $u(x,t)=\frac{1}{\sqrt t}U\Bigl(\frac{x}{\sqrt t}\Bigr)$
is self-similar and the profile $U(x)$ is such that
\begin{equation}
\label{pro2Dss}
U(x)=a(x)-\log(|x|)\FF(x)\!:\!\biggl(\int_{\S^1}a\otimes a\biggr) + \mathcal{O}(|x|^{-3}),
\qquad\hbox{as $|x|\to\infty$}.
\end{equation}
\end{theorem}

\Proof
The first statement follows from the bi-integral formula~\eqref{biff} and our previous Lemmata.
Indeed,
as we have already observed, by Lemma~\ref{as heat lem}, we can write $e^{t\Delta}a(x)=a(x)+\mathcal{O}(t|x|^{-3})$.
Next, writing
$$B(e^{t\Delta}a,e^{t\Delta}a)=\L(e^{t\Delta}a\otimes e^{t\Delta}a),$$
we apply first Lemma~\ref{lle} with $w=e^{t\Delta}a\otimes e^{t\Delta}a$,
and then Lemma~\ref{lol}.
This shows that $-B(e^{t\Delta}a,e^{t\Delta}a)$  equals to the second term on the right hand side of~\eqref{pro2Dnss},
up to an error $\mathcal{O}(t|x|^{-3})$ for large~$|x|$.
The last two terms in the bi-integral formula can also be included into the remainder term
$\mathcal{O}(t|x|^{-3})$,
as shown by combining inequalities~\eqref{qw}-\eqref{qwB} with Lemma~\ref{le5}.

In the case of homogeneous data, an elementary computations shows that
$$\int_0^t\!\!\int_{\sqrt s\le |y|\le |x|} (a\otimes a)(y)\,dy\,ds
=\biggl(\int_{\S^1} a\otimes a\biggr)t\log\Bigl(\frac{|x|}{\sqrt t}\Bigr)+t/2.$$
Then profile~\eqref{pro2Dss} follows from profile~\eqref{pro2Dnss}
passing to self-similar variables and eliminating~$t$.

\endProof

\section{Asymptotics in the higher-dimensional case}
\label{sec9}

We now establish the analogue of Lemma~\ref{lle} for the higher dimensional case.
\begin{lemma}
\label{lem ter}
Let $w=(w_{h,k})$ with $w_{h,k}\in X^1_4$. Then we have, as $|x|\to\infty$, 
and uniformly in time, for $|x|\ge e\sqrt t$,
\begin{subequations}
\begin{equation}
 \label{exp Lw d}
\L(w)(x,t)=\FF(x)\!:\!\int_0^t\!\!\int w(y,s)\,dy\,ds \,+\, \mathcal{O}\bigl(t|x|^{-5}\log(|x|/\sqrt t) \bigr)
\end{equation}
for  $d=3$, and
\begin{equation}
\label{exp Lw dd}
\L(w)(x,t)=\mathcal{O}\bigl(t|x|^{-5}\log(|x|/\sqrt t) \bigr)
\end{equation}
when $d\ge4$.
\end{subequations}
\end{lemma}

\Proof
We go back to the decomposition $\L=\L_1+\L_2+\L_3$ obtained in~\eqref{dec L}.
Writing~$\L_1$ as in~\eqref{L1 dec} and using the estimate $|w(y,s)|\le C(\sqrt s+|x|)^{-4}$,
the bound $|\nabla F(x,t)|\le C|x|^{-d-2}$, 
and the fast decay of $\widetilde \Psi$ shows that the first and the third term
in~\eqref{L1 dec} are bounded by $Ct|x|^{-5}$ (with an additional logarithmic factor $\log(|x|/\sqrt t)$,
for the first term in~\eqref{L1 dec}, when $d=3$) for $|x|\ge e\sqrt t$.
The second term in~\eqref{L1 dec} has the form
$$\FF(x)\!:\!\int_0^t\!\!\int_{|y|\le |x|/2} w(y,s)\,dy\,ds .$$
Using again that $|w(y,s)|\le C(\sqrt s +|x|)^{-4}$ and distinguishing the between the cases $d=3$ and $d\ge4$
shows that such term can be written as the right-hand sides in~\eqref{exp Lw d}-\eqref{exp Lw dd}.

We now decompose  $\L_2$, as
\begin{equation}
\label{L2 decc}
\begin{split}
\L_2 &=  \int_0^t\!\!\int_{|y|\le |x|/2} F(y,t-s)\bigl[w(x-y,s)-w(x,s)\bigr]\,dy\,ds\\
 &\qquad +\int_0^t w(x,s)\int_{|y|\le |x|/2} F(y,t-s)\,dy\,ds. 
\end{split}
\end{equation}
Since $|\nabla w(x,t)|\le C|x|^{-5}$, the first term in~\eqref{L2 decc} is bounded by $C|x|^{-5}t\log(|x|/\sqrt t)$
for $|x|\ge e\sqrt t$. Combining the estimate $|F(y, t-s)|\le |y|^{-d-1}$ with the condition $\int F(\cdot,t-s)\,ds=0$,
shows that the second term in~\eqref{L2 decc} is bounded by~$Ct|x|^{-5}$.
Such bound holds also for $\L_3$ as easily checked using the usual spatial decay estimates of $F$ and $w$.

\endProof

\medskip
Our next Lemma essentially states that if $a$ and $b$ are two functions
defined on~$\R^d$ and well behaved at infinity
(for example, the derivatives of $a$ and $b$ decay faster than $a$ and $b$
as $|x|\to\infty$), then
\begin{equation*}
 \label{fun h}
(e^{t\Delta}a)(e^{t\Delta}b)\sim e^{t\Delta}(ab), \qquad\hbox{as ${|x|\to\infty}$}.
\end{equation*}
More precisely, we have:

\begin{lemma}
\label{fun l}
Let $d\ge3$ and   $a,b\in \dot E^1_1$.
Then
\begin{equation}
 \label{fun rh}
(e^{t\Delta}a)(e^{t\Delta}b)=e^{t\Delta}(ab)-2\int_0^t e^{(t-s)\Delta}\bigl[\nabla e^{s\Delta}a\cdot \nabla e^{s\Delta}b\bigr]\,ds.
\end{equation}
\end{lemma}

\Proof
Let $v=e^{t\Delta}a$ and  $w=e^{t\Delta}b$. Then we have $\partial_t v=\Delta v$ and $\partial_t w=\Delta w$.
Multiplying by $w$ the first equation and by $v$ the second one we get
\begin{equation*}
 \partial_t(vw)=w\Delta v+v\Delta w=\Delta(vw)-2\nabla v\cdot\nabla w.
\end{equation*}
Since $d\ge3$, $ab$ is locally integrable in $\R^d$.
But $(vw)(t)\to ab$ as $t\to0^+$ weakly (because $v(t)\to a$ and $w(t)\to b$ in $L^2_{\rm loc}(\R^d)$,
for example, as $t\to0$).
Then the conclusion follows from Duhamel formula.

\nobreak
\endProof

\medskip
In the above Lemma we only used, in fact, $a,b\in E^0_1$. The stronger assumption
$a,b\in \dot E^1_1$, however, ensures that the last term in~\eqref{fun rh},
decays faster as $|x|\to\infty$ than $e^{t\Delta}(ab)$.

\medskip
We now give the higher-dimensional counterpart of Theorem~\ref{theo2D}.

\begin{theorem}
\label{theorem3d}
Let $u(x,t)\in X^3_1$ be the global solution of the Navier-Stokes equations starting from
$a\in \dot E^3_1$ (as constructed in Proposition~\ref{theorem1}).
Then $u$ has the following profile as $|x|\to\infty$, uniformly in time for  $|x|\ge e \sqrt t$.
For $d=3$,
\begin{subequations}
\begin{equation}
 \label{pro3Dnss}
u(x,t) =e^{t\Delta}a(x)-t\,e^{t\Delta}\P\nabla\cdot(a\otimes a)-\FF(x)\!:\!\Lambda(t) 
 +\mathcal{O}\biggl(t^2|x|^{-5}\log\Bigl(\frac{|x|}{\sqrt t}\Bigr)\biggr),
\end{equation}
for some matrix-valued function $\Lambda(t)=(\Lambda_{h,k}(t))$,
satisfying  $|\Lambda(t)|\le Ct^{3/2}$.
Moreover, when $d\ge4$,
\begin{equation}
\label{pro4Dnss}
u(x,t) = e^{t\Delta}a(x)-t\,e^{t\Delta}\P\nabla\cdot(a\otimes a) 
 +\mathcal{O}\biggl(t^2|x|^{-5}\log\Bigl(\frac{|x|}{\sqrt t}\Bigr)\biggr).
\end{equation}
\end{subequations}
\end{theorem}

\begin{remark}
The function $\Lambda(t)$ is not know explicitly, but it depends on $u$ and $a$ in an explicit way:
see formula~\eqref{ssos} below.
For more regular data, namely $a\in \dot E^4_1$, and recalling Lemma~\ref{as heat lem}
 (applied with $m=4$ and $\vartheta=1$) one can replace in the above asymptotics
the term $e^{t\Delta}a(x)$ with 
$a(x)+t\Delta a(x)$.

\end{remark}

\Proof
As for the proof of our previous theorem, we write $u$
by means of the bi-integral formula~\eqref{biff}.
As an application of Lemma~\ref{fun l} we can rewrite (for $d\ge3$)
the term $B(e^{t\Delta}a,e^{t\Delta}a)$ appearing
in the bi-integral formula~\eqref{biff} in a more convenient form
(we denote here by ${}^T\!\!A$ the transposed of the matrix $A$):
\begin{equation}
 \label{Baa}
\begin{split}
B(e^{t\Delta} &a  ,e^{t\Delta}a) = \int_0^t F(t-s)* \bigl(e^{s\Delta}a\otimes e^{s\Delta}a\bigr)\,ds\\
&=\int_0^t e^{(t-s)\Delta}\P e^{s\Delta}\nabla\cdot(a\otimes a)
-2\int_0^t\!\!\int_0^s e^{(t-\tau)\Delta}\P\nabla\cdot\Bigl[\!{\phantom{\Bigl|}}^T\!\! \Bigl(\nabla\otimes e^{\tau\Delta}a\Bigr)
	\Bigl(\nabla\otimes e^{\tau\Delta}a\Bigr)\Bigr]\,d\tau\,ds\\
&= t e^{t\Delta}\P\cdot\nabla(a\otimes a)-2\int_0^t(t-\tau) e^{(t-\tau)\Delta}\P\nabla\cdot
   \Bigl[ \!{\phantom{\Bigl|}}^T\!\! \Bigl(\nabla\otimes e^{\tau\Delta}a\Bigr) \Bigl(\nabla\otimes e^{\tau\Delta}a\Bigr) \Bigr]\,d\tau\,,
\end{split}
\end{equation}
where we applied Fubini's theorem in the last equality.

\medskip

We set
$$
\widetilde\L(w)(t)\equiv\int_0^t(t-\tau)F(t-\tau)*w(\tau)\,d\tau
 \quad\hbox{and}\quad
\overline\L(w)(t)\equiv\int_0^t \tau\,F(t-\tau)*w(\tau)\,d\tau
  .$$
Note that, excepted for the additional factors $t-\tau$ or $\tau$, the operator $\widetilde{\L}$
and $\overline{\L}$ agree with the operator
$\L$ introduced in~\eqref{Lw} and studied before.
If we introduce the matrix
$$
    w_1\equiv 
    \!{\phantom{\Bigl|}}^T\!\! \Bigl(\nabla\otimes e^{\tau\Delta}a\Bigr) \Bigl(\nabla\otimes e^{\tau\Delta}a\Bigr),
$$
then we can rewrite~\eqref{Baa} as
\begin{equation*}
 B(e^{t\Delta}a,e^{t\Delta}a)=t\,e^{t\Delta}\P\nabla\cdot(a\otimes a)-2\widetilde\L(w_1).
\end{equation*}
The estimates of Lemma~\ref{heat l} (in the case $m=\vartheta=1$), imply $w_1\in X^1_4$.
But the result of Lemma~\ref{lem ter}, established before for the operator $\L$, can be easily adapted
to the operators $\widetilde\L$ and $\overline\L$; indeed the factors $t-\tau$ and $\tau$ are 
harmless in our estimates due to the obvious inequalities $t-\tau\le t$ and $\tau\le t$.
Thus, we get, for $d=3$,
\begin{equation*}
 \widetilde\L(w_1)=
\FF(x)\!:\!\int_0^t\! (t-\tau)
 \!\!\int  w_1 \,dy\,ds
 \,+\,\mathcal{O}\Bigl( t^2|x|^{-5}\log\bigl(|x|/\sqrt t\bigr)\Bigr).
\end{equation*}

When $d\ge4$, we can simply write
\begin{equation*}
 \widetilde\L(w_1)=\mathcal{O}\bigl(t^2|x|^{-5}\log(|x|/\sqrt t)\bigr).
\end{equation*}

\medskip
It remains to write the asymptotics (or to estimate) the two last terms 
$B(e^{t\Delta}a,B(u,u))$ and
$B(B(u,u),B(u,u))$ appearing
in the bi-integral formula~\eqref{biff}.
Let 
$$w_2\equiv\tfrac{1}{t} \,e^{t\Delta}a\otimes B(u,u).$$
We get from Lemma~\ref{heat l} (applied with $m=1$ and $\vartheta=1$)
and Remark~\ref{rem fine} that $w_2\in X^1_4$.
In the same way,
Remark~\ref{rem fine} ensures that, if we set 
$$w_3\equiv \tfrac{1}{t} \,B(u,u)\otimes B(u,u),$$
then  $w_3\in X^1_4$.
Therefore, Lemma~\ref{lem ter} (or more precisely, the adaptation of this Lemma to $\overline\L(w_2)$
and $\overline\L(w_3)\,$)  implies, for  $d=3$,
\begin{equation}
\label{babuu}
\begin{split}
2B(e^{s\Delta}a, & B(u,u))-B(B(u,u),B(u,u)) \\
&=\,\,\,2\overline\L(w_2)-\overline\L(w_3)\\
&=\,\,\, \FF(x)\!:\!\int_0^t s\!\!\int (2w_2-w_3)\,dy\,ds +\mathcal{O}\Bigl(t^2|x|^{-5}\log(|x|/\sqrt t)\Bigr),\\
\end{split}
\end{equation}
as $|x|\to\infty$.

\medskip
When $d\ge4$, the first term in the right-hand side of~\eqref{babuu}
can be dropped.
Therefore,
the proof of the expansion~\eqref{pro4Dnss} follows from the bi-integral formula,
collecting the above estimates.

\medskip
In the case $d=3$, it is convenient to introduce the time-dependent matrix
\begin{equation}
 \label{ssos}
\Lambda(t)=\int_0^t\!\!\int \Bigl[ -2(t-s)w_1-2s\, w_2+s\, w_3\Bigr]\,dy\,ds.
\end{equation}

The expansion~\eqref{pro3Dnss} now follows by collecting all the above expressions.
The estimate $|\Lambda(t)|\le Ct^{3/2}$ is immediate,
because $w_1$, $w_2$ and $w_3$ belong to $X^1_4$.

\endProof

\medskip

As an application of this theorem, we can complete the proof of Theorem~\ref{theoss}
by giving the far-field asymptotics of self-similar solutions in the case $d\ge3$.

\bigskip
\noindent
{\it End of the Proof of Theorem~\ref{theoss}.\/}
We assumed that $a\in C^\infty(\S^{d-1})$ and that $a$ is is homogeneous of degree~$-1$.
From the second part of Lemma~\ref{as heat lem},
$$ e^{t\Delta}a(x)=a(x)+t\Delta a(x)+\mathcal{O}(t^2|x|^{-5}).$$
But the solution $u$ is of the self-similar form $u(x,t)=\frac{1}{\sqrt t}U(x/\sqrt t)$.
Moreover, the linear part $e^{t\Delta}a$ and the nonlinear part $B(u,u)$ of $u$
are also of self-similar form, so that, with the same notations 
of the previous proof, $w_j(y,s)=\frac{1}{s^2}W_j(y/\sqrt s)$, where
$$W_j(y)=w_j(y,1), \qquad j=1,2,3.$$

If follows from~\eqref{ssos} that, in the case $d=3$,  $\Lambda(t)$ is of the form $\Lambda(t)=t^{3/2}B$,
for some \emph{constant\/} matrix $B=(B_{h,k})$.
As for $\Lambda(t)$, such matrix $B$ is not known explicitly, however, it is possible
to obtain an explicit integral formula
relating $B$ to the datum $a$ and the profile $U$,
performing a self-similar change of variables in the integral~\eqref{ssos}.
An easy computation yields
\begin{equation}
\label{Bmat}
B=\frac{1}{3}\int \Bigl(-8W_1-4 W_2+2 W_3\Bigr)(y)\,dy.
\end{equation}

\begin{subequations}
Now we can pass to self-similar variables in expansion~\eqref{pro3Dnss} and, after eliminating~$t$,
we get, for  $d=3$,
\begin{equation}
\label{pro3Dp}
U(x)=a(x)+\Delta a(x)- e^{\Delta}\P \cdot\nabla (a\otimes a)
-\frac{ \mathcal{Q}(x)\!:\!B}{|x|^7} +\mathcal{O}\bigl(|x|^{-5}\log(|x|)\bigr),
\end{equation}
 as $|x|\to\infty$.

As before, for $d\ge4$, the far-field asymptotics has a simpler structure, namely,
\begin{equation}
\label{proDp}
U(x)=a(x)+\Delta a(x)- e^\Delta\P \cdot\nabla (a\otimes a)+ \mathcal{O}\bigl(|x|^{-5}\log(|x|)\bigr),
\end{equation}
as $|x|\to\infty$.
\end{subequations}

\medskip
To finish the proof, it remains to show that we can drop the filtering operator $e^\Delta$ 
appearing in the right-hand side of equations~\eqref{pro3Dp} and~\eqref{proDp}.
Recall that $a$ is smooth on the sphere. In fact, the condition $a\in C^\infty(\S^{d-1})$
will allow us to carry the proof using only ``soft arguments''.
The datum $a$ being homogeneous of degree~$-1$,
 $\nabla\cdot (a\otimes a)$ is a homogeneous distribution of degree~$-3$
(here we need $d\ge3$), which
agree with a $C^\infty$ function outside the origin.
But the matrix Fourier multiplier of the operator $\P$ (given by $\delta_{j,k}-\xi_j\xi_k|\xi|^{-2}$)
is homogeneous of degree zero and smooth outside the origin).
Then it follows (see, {\it e.g.\/}, \cite[p. 262]{Ste93})
that $\P\nabla\cdot(a\otimes a)$ is a homogeneous distribution of degree~$-3$
that agrees with a $C^\infty$ function outside the origin.

Now let $\chi\in C^\infty_0(\R^d)$ be a cut-off function equal to~$1$ in a neighborhood of the origin
and write 
$$e^{\Delta}\P\nabla\cdot(a\otimes a)=e^{\Delta}\chi \P\nabla\cdot(a\otimes a)+
  e^{\Delta}(1-\chi)\mathcal{A}(x),$$
where $\mathcal{A}(x)$ a smooth function on $\R^d$, agreeing with $\P\nabla\cdot(a\otimes a)$
outside a neighborhood of the origin.
In particular, $(1-\chi)\mathcal A\in E^m_3$, for all $m\in\N$.

Note that
$e^{\Delta}\chi\nabla\cdot(a\otimes a)$ is an analytic function, given by
$$ e^{\Delta}\chi\P\nabla\cdot(a\otimes a)(x)=\langle\chi\P\nabla\cdot(a\otimes a),g_1(x-\cdot)\rangle,$$
where $g_1$ the standard gaussian and the $\langle\cdot,\cdot\rangle$
refers to the duality product between
compactly supported distributions and 
$C^\infty$ functions.
The properties of compactly supported distributions guarantee
the existence of a compact $K$ in $\R^d$ and  $C>0$, $M\in\N$ such that
$$ \Bigl| \langle\chi\P\nabla\cdot(a\otimes a),g_1(x-\cdot)\rangle \Bigr|
    \le C\sum_{|\alpha|\le M}\sup_{y\in K}\partial^\alpha g_1(x-y)\le C'g_1(x/2)$$
for large enough~$|x|$.
In particular,
$e^{\Delta}\chi\nabla\cdot(a\otimes a)=\mathcal{O}(|x|^{-5})$ as $|x|\to\infty$.

Let us now apply the asymptotic formula for convolution integrals~\eqref{cpp}
with $g=g_1$ and $f=(1-\chi)\mathcal{A}$.
We obtained this formula under the assumption $f\in \dot E^m_\vartheta$, with
$0\le\vartheta<d$. Here we have, instead, $f\in E^m_3\subset \dot E^m_3$ but it is
easily checked that such formula remains valid, in this case,  also when $d=3$, with the same proof,
since $f$ is locally integrable.
Applying this formula in the case $m=2$, and using $\int g_1=1$ and $\int y\,g_1(y)\,dy=0$,
we get, for $|x|\to\infty$,
\begin{equation*}
\begin{split}
 e^{\Delta}(1-\chi)\mathcal{A}(x)&=g_1*f(x)=f(x)+\mathcal{O}(|x|^{-5})\\
 &=\mathcal A(x) +\mathcal{O}(|x|^{-5})\\
   & =\P\nabla\cdot(a\otimes a)(x)+\mathcal{O}(|x|^{-5}).
\end{split}
\end{equation*}
Theorem~\ref{theoss} is now completely proved.

\endProof


\section{Acknowledgements}
The preparation of this paper
was supported in part by the European Commission Marie Curie Host Fellowship
for the Transfer of Knowledge ``Harmonic Analysis, Nonlinear
Analysis and Probability''  MTKD-CT-2004-013389,
and in part by the program Hubert-Curien ``Star'' N.~16560RK.

\bibliographystyle{amsplain}

\end{document}